\documentstyle{amsppt}
\document

\topmatter
\title
HyperK\"ahler Manifolds and Birational Transformations in dimension 4
\endtitle
\author Dan Burns, Yi Hu and Tie Luo
\endauthor
\date 
\enddate
\subjclass 
\endsubjclass
\address
 Department of Mathematics, University of Texas, Arlington, TX 76019
\endaddress
\email
hu\@math.uta.edu, tluo\@math.uta.edu
\endemail
\abstract
We show that any birational map between projective
hyperK\"ahler manifolds of dimension 4 is composed
of a sequence of simple flops or elementary Mukai 
transformations under the assumption that each irreducible 
component of the indeterminacy of the birational map is normal.
\endabstract

\endtopmatter
\openup6pt

\heading
\S 1 Introduction and statement of main theorem
\endheading

The main result of this paper is a solution to an open problem posed by Mukai 
more than a decade ago
(Problem 4.5, [Mu2]) under a normality assumption.

\proclaim{Theorem 1.1}
Let $\Phi: X --\rightarrow X'$ be a birational transformation 
between two nonsingular hyperK\"ahler  
fourfolds and $B \subset X$ the indeterminacy of $\Phi$. Assume each irreducible component of $B$ is normal. Then $B$ 
 is a union of a $\Bbb P^2$ and some rational surfaces which are either $\Bbb P^2$s or blowups of  $\Bbb P^2$s and  $\Phi: X --\rightarrow X'$ can be decomposed as a 
sequence of the Mukai elementary transformations along these $\Bbb P^2$s up to an isomorphism.
\endproclaim

This theorem gives a complete classification of birational transformations of  
projective symplectic  fourfolds. Note that the birational maps between Calabi-Yau
threefolds have been classified in [Ko1].

A hyperK\"ahler manifold is a projective manifold $X$
equipped with a holomorphic symplectic form $\omega$. Such a manifold has trivial canonical bundle.  
It is desired, following Mori's minimal model program,
that any birational map between two  minimal models with trivial canonical bundles
can be decomposed as a finite sequences of elementary ones, flops.
A particularly interesting 
class of flops consists of the so-called Mukai's elementary transformations [Mu2].

Precisely, let $(X, \omega)$ be a holomorphic symplectic manifold  and $P$ an embedded
smooth subvariety. Assume further that $P$ is a $\Bbb P^r$-bundle over a smooth
variety $\Sigma$ such that the codimension of  $P$ coincides with $r$. Then it is known [Mu2]
that one can blow up $X$ along $P$ to get a smooth variety $\widetilde X$ and the exceptional
divisor can be blow down along a different ruling to get a smooth variety $X'$.
Moreover, $X'$ comes equipped with a symplectic form $\omega'$ which coincides
with $\omega$ away from the exceptional locus. 
Such a simple birational process $X --\rightarrow X'$ is 
called a Mukai elementary transformation.

The first nontrivial dimension that such elementary transformations can occur is 4.
Here, $P$ is necessarily the projective space $\Bbb P^2$.
Trying to find a solution to Mukai's open problem,
we wondered whether an irreducible flop is necessarily a
Mukai's elementary transformation. This turns out to be true in the projective category assuming normality of the indeterminacy (Theorem 1.1).

It is clear that one needs to isolate a copy of $\Bbb P^2$ to be able to perform a Mukai elementary transformation. For this we need to classify the exceptional locus of a logterminal contraction, which contracts part of the indeterminacy. 
It turns out to be a formidable task to fully achieve this goal and many years have past since we started working on
the project. Even if assuming the normality of each irreducible component, the task is still quite involved.
The relative relation of various components can be very complex.
One key point is to check that the normality of each irreducible component in $B$ is kept after METs.
To this end, a useful result of Wierzba [W] is applied.

There are examples of birational maps between symplectic fourfolds with indeterminacy consisting a chain of rational surfaces as shown below.

\demo{Example} Let $S$ be the K3 surface defined by the quartic homogeneous polynomial

$$
xy(x^2+y^2+z^2+w^2)+w(x^3+y^3+z^3+w^3)
$$
in $\Bbb P^3$ with coordinates $[x,y,z,w]$. Note that there are three $\Bbb P^1$s on $S$ in the $\Bbb P^2$ defined by $w=0$, i.e; two lines and a conic.
Let $S^{[2]}$ be the Hilbert scheme of points of length $2$ on $S$. There is a chain of three $\Bbb P^2$, say $B_1, B_2, B_3$ on $S^{[2]}$ corresponding to the three ${\Bbb P}^1$ on $S$. Let $B_1, B_3$ be the $\Bbb P^2$s corresponding to the two lines. $B_1, B_3$ are disjoint and they both have one point in common with $B_2$. Now we perform a MET to $B_2$ to get $Y$ and METs to $B_1, B_3$ to get $X$. Let $f$ be the birational map from $X$ to $Y$. The indeterminacy on $X$ is a chain of three surfaces with the two ends isomorphic to $\Bbb P^2$s and the middle one isomorphic to blowup of $\Bbb P^2$ at two points. The indeterminacy on $Y$ is a chain of three surfaces with two ends isomorphic to blowup of $\Bbb P^2$ at a point and the middle one isomorphic to $\Bbb P^2$.
\enddemo
 
The example indicates that our situation is more complicated than the picture demonstrated in [Ka1]
where the indeterminacy is a disjoint union of $\Bbb P^2$s. But,
it is easy to prove that any two components that are isomorphic to $\Bbb P^2$ are either disjoint
or meeting at isolated points. For otherwise, they would meet in a locus of dimension 1. Contracting
one of the component in $X$ (being isomorphic to $\Bbb P^2$, it is contractible) will result in
the contraction of the dimension 1 locus in the other component
(also $\cong \Bbb P^2$), which is absurd.

\demo{Acknowledgments}. We thank 
J\'anos Koll\'ar for his help,  Nick Shepherd-Barron for pointing out an error,
Jan Wierzba for sharing his results in [W]. We also thank  S.-T. Yau for being
interested in this work.

\enddemo

\heading
\S 2 Some General Lemmas
\endheading

Throughout the paper $X$ stands for a symplectic fourfold unless otherwise stated.

\proclaim{Definition 2.1}
Let $X$ be a symplectic fourfold. A MET from
$(X,B)$ to $(X',B')$ is the Mukai elementary transformation
$$
\matrix
 & & E\subset \bar X& &  \\
& &p\swarrow\quad   \searrow q&  & \\
&B\simeq \Bbb P^2\subset X &--\rightarrow & X'\supset B'\simeq (\Bbb P^2)^*&
\endmatrix
$$
where $E$ is the incidence correspondence between $B$ and $B'$.
\endproclaim

Regarding a MET, one has the following basic result.

\proclaim{Lemma 2.2}
Let $f: X --\rightarrow X'$ be a MET with exceptional loci $B\subset X$,
$B'\subset X'$. Let $E\subset \bar X$ be the exceptional divisor. Let
$H$ be an divisor such that $H\cdot C<0$ for any curve $C\subset B$ and $H'$ the
proper transform of $H$ in $X'$. Then $p^*H-q^*H'=aE$ for $a>0$. Moreover, $H'$ is
numerically positive on $B'$, i.e., $H'\cdot C'>0$ for any curve $C'\subset B'$.
\endproclaim

\demo{Proof }
Since $E$ is the only exceptional divisor for both $p$ and $q$, we have
$p^*H-q^*H'=aE$. To see the sign of $a$, let $C$ be a line of $B$. Recall
that $B'$ is the dual space of $B$. We set $\bar C=(C,P)\subset E$ 
 where $P \in B'$ is the point that corresponds to $C$. Thus we have
$$
aE\cdot \bar C=p^*H\cdot \bar C-q^*H'\cdot \bar C=H\cdot p_*\bar C-H'\cdot q_*\bar C=H\cdot C<0
$$
by the projection formula. This implies $a>0$ since $E\cdot \bar C<0$.

To get the last statement,
it suffices to show that $H'\cdot C'>0$ for any line $C'$ in $B'$. Like above,
set $\bar {C'}=(Q',C')\subset E$ where $Q'$ is the point in $B$ that corresponds to
$C'$. Then,
$$
0>aE\cdot \bar {C'}=p^*H\cdot \bar {C'}-q^*H'\cdot \bar {C'}=H\cdot p_*\bar {C'}-H'\cdot q_*\bar
{C'}=-H'\cdot C'.
$$
\qed
\enddemo

Several technical results are also needed. 

First, one expects that the exceptional locus of $\Phi: X --\rightarrow X'$ contains a rational curve.
The following lemma asserts that in the case of symplectic variety, any rational curve moves in 
one more family than Riemann-Roch predicts.

\proclaim{Lemma 2.3 (Ran)} 
 Assume a symplectic manifold $X^n$ contains a rational curve $C$, 
then $C$
deforms at least in $(n-2)$ families.
\endproclaim

\demo{Proof} See the proof in [R] when $C$ is smooth.

When $C$ is singular, we consider the graph of $f:\Bbb P^1\rightarrow
C\subset X$:
$$
\bar f: \Bbb P^1\rightarrow \bar C\subset \Bbb P^1\times X=\bar X.
$$
$\bar C$ is smooth. Let $\bar{\Cal H}$ ($\Cal H$) be the Hilbert scheme
containing $\bar C$ in $\bar X$ ($C$ in $X$), one has the following
estimate by [R] (see also [Ka2])
$$
\text{dim}\bar{\Cal H}\geq \chi (N_{\bar C/\bar X})+\text{dim im}(\pi) 
$$
$$
= -\bar C.K_{\bar X}+(n+1-3)+\text{dim im}(\pi) $$

$$= 2+n+1-3+\text{dim
im}(\pi)$$
$$ = n+\text{dim im}(\pi)
$$
where $\pi$ is the semi-regularity map, whose dual is
$$
\pi^t: H^0(\bar X,\Omega^2_{\bar X})\rightarrow H^0(\bar C,N^*_{\bar
C/\bar X}\otimes\Omega_{\bar C}).
$$
To see that $\pi^t$ is nontrivial, we consider the image of
$\beta^*\omega$ where $\omega$ is the holomorphic symplectic form on $X$
and $\alpha, \beta$ are projections from $\bar X$ to $\Bbb P^1$ and $X$.
We shall show that $\beta^*\omega$ is not zero at any point
$y=(x,f(x))\in \bar C$ as long as $f(x)$ is a smooth point on $C$.
Around an analytic neighborhood of $y$ in $\bar X$ which is viewed as the
product of analytic neighborhoods of $x$ and $f(x)$ in $\Bbb P^1$ and
$X$ (respectively), one has a (non-canonical) isomorphism
$$
\Omega^2_{\bar X}\simeq
\beta^*\Omega^2_X+\beta^*\Omega_X\wedge\alpha^*\Omega_{\Bbb P^1}
$$
$$\simeq \beta^*\wedge^2 N^*_{C/X}+\beta^*(\Omega_{C}\otimes N^*_{C/X})
+\beta^*N^*_{C/X}\otimes\alpha^*\Omega_{{\Bbb
P}^1}+\beta^*\Omega_C\wedge\alpha^*\Omega_{\Bbb P^1}.
$$
Under this identification, we have
$$
\wedge^2N^*_{\bar C/\bar X}\simeq
\beta^*\wedge^2N^*_{C/X}+\beta^*N^*_{C/X}\otimes\alpha^*\Omega_{{\Bbb
P}^1}.
$$
Then the non-degenerate property of $\omega$ implies that the component of
$\beta^*\omega$ in $\beta^*(N^*_{C/X}\otimes\Omega_C)$ is not trivial.
Thus
$\pi^t(\beta^*(\omega))\neq 0$.

So the semi-regularity map $\pi$ is nontrivial. Hence 
$$
\text{dim}\bar{ \Cal H }\geq n+1.
$$
Notice that a nontrivial automorphism of $\Bbb P^1$ from the first
factor of $\bar X$ gives a nontrivial deformation of $\bar C$, which however
does not move $C$ in $X$. Therefore
$$
\text{dim}{\Cal H}\geq n+1-3=n-2.
$$
\qed
\enddemo

When $n=4$, we obtain that each and every rational curve in a symplectic
fourfold moves in a at least 2-dimensional family.

The next lemma was pointed out by J. Koll\'ar.

\proclaim{Lemma 2.4}Let $S$ be a normal surface, proper over $\Bbb C$.
Then $S$ satisfies exactly one of the following:
\roster
\item  Every morphism $f: {\Bbb P}^1 \rightarrow S$ is constant;
\item  There is a morphism $f: {\Bbb P}^1 \rightarrow S$ such that
       $f$ is rigid.
\item  $S \cong {\Bbb P}^2$;
\item  $S \cong {\Bbb P}^1 \times {\Bbb P}^1$, or $S$ is isomorphic to
       a minimal ruled surface over a curve of positive genus
       or a minimal ruled surface with a negative section
       contracted.
\endroster
\endproclaim

\demo{Proof } If we have either of (1) or (2), we are done.
Otherwise, there is a morphism $f: \Bbb P^1 \rightarrow S$ deforms
in a 1-parameter family, thus $S$ is uniruled.

Let $p: \bar S \rightarrow S$ be the minimal desingularization
with the exceptional curve $E$. $\bar{S}$ is also uniruled,
hence there is an extremal ray $R$. There are 3 possibilities for $R$.

\roster
\item  $\bar S \cong \Bbb P^2$, thus also $S \cong \Bbb P^2$ (which implies (3)).
\item  $\bar{S}$ is a minimal ruled surface (which implies (4)).
\item  $R$ is spanned by a (-1)-curve $C_0$ in $\bar S$.
\endroster

But (3) is impossible, because the image of $C_0$ in $S$ would have
been rigid. The proof goes as follows. Assume the contrary that
$f_0: \Bbb  P^1 \cong C_0 \subset \bar S \rightarrow S$ is not rigid
and let $f_t: \Bbb P^1 \rightarrow S$ be a 1-parameter deformation.
For general $t$, $f_t$ lifts to a family of morphisms
$\bar{f}_t : \Bbb P^1 \rightarrow \bar S$. As $t \to 0$,
the curves $\bar{f}_t (\Bbb P^1)$ degenerate and we obtain a cycle
$$\lim_{t \to 0} \bar{f}_t (\Bbb P^1) = C_0 + F$$
where $\text{ Supp} F \subset \text{ Supp} E$.
$\bar{S}$ is the minimal resolution, thus $K_{\bar{S}} \cdot F \ge 0$.
Therefore,
$$K_{\bar{S}} \cdot \bar{f}_t ({\Bbb P}^1) \ge K_{\bar{S}} \cdot C_0 = -1.$$
On the other hand, for a general $t$ the morphism $\bar{f}_t$ is free,
thus
$$K_{\bar{S}} \cdot \bar{f}_t ({\Bbb P}^1) \le -2$$
by II.3.13.1, [Ko2]. This contradication shows that
$f_0$ is rigid.
\qed
\enddemo

We will also use
the following lemma which is essentially from 2.19 of [Ketal].
 
\proclaim{Lemma 2.5} 
Let $\Phi :X--\rightarrow X'$ be a birational map between projective symplectic fourfolds.
Assume $H'$ is ample on $X'$ and $H$ its
proper transform on $X$. $\Phi$ is a morphism if $H$ is nef. $\Phi$ is an
isomorphism if $H$ is ample (or numerically positive).
\endproclaim

\demo{Proof} See the proof of Proposition 2.7 of [C] for the details when $H$ is ample.
The same proof goes through when $H$ is numerically positive.
\qed
\enddemo

Finally, the following proposition due to J.~Wierzba [W] will be very useful for checking the invariance of the normality of the exceptional locus under METs.

\proclaim{Proposition 2.6} (Theorem 1.3 and 1.4 in[W])
Let $\pi:\hat X\rightarrow X$ be an isolated symplectic singularity of dimension $4$ and $E$ the exceptional locus. $E$ is a union of irreducible projective surfaces $B_i$ whose normalizations are $\Bbb P^2$. Then

1) if $E_i$ meets $E_j$ along a curve $C$, $C$ is in the singular locus of either $E_i$ or $E_j$;

2) if $E_i$ is nonsingular in codimension $1$ for all $i$, $E_i$ is normal for all $i$.

\endproclaim

\heading
\S 3 Proof of the main theorem
\endheading

\demo{Proof} We now start to prove Theorem 1.1. To begin with,
let $$\Phi: X--\rightarrow X'$$ be a birational map between two hyperK\"ahler fourfolds.
This map is necessarily isomorphic in codimension 1.  
Our goal is to show that the indeterminancy of $\Phi$  is a
union of a ${\Bbb P}^2$ with other rational surfaces and $\Phi$ is factored into METs.

Let $H'$ be a very ample divisor on $X'$. Let $H$ be its
proper transform in $X$.

We divide the proof into a few steps.

1. First, we consider the pair $(X,\epsilon H)$. It is log-terminal for
$\epsilon<<1$. Since $H$ is not nef on $X$
(actually not nef on $B=\cup B_i$, the union of all irreducible subvarieties where
$f$ is not defined), by Lemma 2.5, there is a curve $C\subset B$ such that $C \cdot H<0$.
Using the contraction theorem [KMM] to the 
log-terminal pair $(X,\epsilon H)$,
there is a morphism $g: X \rightarrow Y$ whose exceptional locus is contained in $B$.
Next, apply the rationality theorem of Kawamata [Ka2] to the morphism $g$, 
the exceptional locus of $g$ is covered by rational curves. 
By abusing notation a little, we still use the letter $C$ to denote a rational curve contracted by $g$. By Lemma 2.2 $C$ moves in at least two famlies and it can not move out of $B$ because $C\cdot H<0$.
Let $B_1$ be an irreducible component of $B$ which is generically swept out by $C$. 
$B_1$ is contracted by the map $g$ that  contracts $C$.
We argue in the following lemma that $B_1$ has to be contracted to a point.

\proclaim{Lemma 3.1} $B_1$ is contracted to a point by the map $g$. (In particular, the resulting
variety $Y$ has only  an isolated singular point.)
\endproclaim

\demo{Proof} Here the proof uses holomorphic Hamiltonian flows. 

First note that by Kawamata (Theorem 2, [Ka2]),
the exceptional locus of $g: X \rightarrow Y$ is covered by a families of 
rational curves. $B_1$ is actually covered by at least a two dimensional family of rational curves. Hence $B_1$ is unirational and thus rational. Clearly, $B_1$ is
a surface. This implies that
$B_1$ is (generically) Lagrangean.

Assume the contrary that the map $g$  mapped $B_1$ onto a curve $D$ in $Y$, 
rather than onto a point as we wish.  Let $f$ 
be a holomorphic function defined in a neighborhood of a general point of $D$, and 
which has $df \ne 0$ when restricted to the curve $D$ (locally around the point). We pull 
this function back up to a neighborhood of a fiber $F$ in the original variety $X$. 
Let $H_f$ be the Hamiltonian holomorphic vector field determined by $f$ and the symplectic 
structure in a neighborhood of the fiber curve $F$. Since the differential 
of $f$ on $B_1$ near $F$ is nonzero, by what we assumed about $f$ and $D$, and since 
$B_1$ is Lagrangean, 
it follows that $H_f$ is transverse to $B_1$ along $F$. Now flowing 
$F$ along the integral curves of $H_f$, we will get a holomorphic deformation of $F$ 
outside of $B_1$, which contradicts the fact that the contraction of the 
extremal curve $F$ contracted only $B_1$ locally around $F$.
\qed of the lemma.
\enddemo

It could happen that a rational curve in $B_1$ moves into another 
component which is also contracted by $g$.
 To show the normalization of $B_1$ is $\Bbb P^2$, 
we need to know additionally that every rational curve moves within $B_1$.  Assume otherwise.
Let $C'$ be a rigid rational curve in $B_1$ and $C''$ a general 
rational curve in $B_1$ meeting $C'$. Then $C''$ moves in only one family which is against Lemma 2.3. 
Hence every rational curves moves in $B_1$.
By Lemma 2.4, there is a morphism 
$$\nu:\Bbb P^2 \rightarrow B_1$$ 
such that $\nu$ is the normalization of $B_1$. So $B_1\simeq \Bbb P^2$.

We now continue from the step 1 of the proof of the main theorem.

2. Next, we perform a MET to $(X,B_1)$ to get $(X^1,B_1^1)$. Let $H^1$ be the
proper transform of $H$. By Lemma 2.2, $H^1$ is numerically positive on
$B_1^1$. Let $\cup_{i \geq 1}B_i^1$ be the image of $\cup_{i\geq 1}B_i$
under the MET. We are done if $H^1$ is numerically positive. Otherwise there is a curve $C^1\subset \cup_{i\geq 1}B_i^1$
such that $C^1 \cdot H^1<0$. Obviously $C^1$ is not contained in $B^1_1$ because of
the positivity of $H^1$ on $B^1_1$. Again, we
apply the contraction-rationality theorem of  [Ka2]
to the log-terminal pair $(X^1,\epsilon_1H^1)$
for $\epsilon_1<<1$ to get a rational curve $C_1$ which deforms in an at least two-dimensional
family. Let $B_2^1$ ($\ne B_1^1$) be the irreducible component which contains the
family. 

$B_2^1$ normalizes to a ${\Bbb P}^2$ and is a ${\Bbb P}^2$ if it is normal by the same arguement used before. There are two situations: 1)
 $B_2^1$ does not intersect $B_1^1$. 2) $B_2^1$ does intersect $B_1^1$. In case 1) $B_2\simeq B_2^1$, so $B_2^1$ is normal.

In case 2), the intersection must be  a set of
finitely many points since $H^1$ is positive on $B_1^1$ and negative on
$B_2^1$. This implies that $B_1$ intersects $B_2$ along some ${\Bbb
P}^1$s before the MET. Moreover the nonnormal locus of $B_2^1$ is contained in $B_1^1\cap B_2^1 $, hence isolated.  If $B_i^1$ is another irreducible component which is also contracted, again $B_i^1$ intersects $B_1^1$ at isolated points. This implies that $B_i^1$ could be nonnormal only at finitely many points since the original $B_i$ is assumed to be normal. Proposition 2.6 says that nonnonormal locus must be empty. After all $B_2^1\simeq {\Bbb P}^2$.

3. Perform a MET to $(X^1,B_2^1)$ to get $(X^2, B_2^2)$. After $k$ steps of doing MET, we arrive at $(X^k, B_k^k)$. If the proper transform $H^k$ of $H$ is positive on $B^k$, we are done. Otherwise there is some $B^k_{k+1}$ whose normalization is ${\Bbb P}^2$ and $H^k$ is negative on it. We may assume that $B^k_{k+1}$ is contracted by a logterminal contraction assoiated with $(X^k,\epsilon_k H^k)$ for $\epsilon_k<<1$. A comparison between $B_{k+1}\subset X$ (which is assumed to be normal) and $B^k_{k+1}$ shows that the nonnormal locus of $B^k_{k+1}$ is a set of finitely many points. More precisely, assume $B^k_{k+1}$ is singular along a curve $C^k$, we backtrack to a previous $i$-th step ($i>1$) after which $B^i_{k+1}$ become singular along the curve $C^i$, the proper transform of $C^k$. Note that $H^i.C<0$ and $B^i_{k+1}$ has to intersect $B^i_i$. Since $H^i$ is positive on $B_i^i$, we conclude that $B^i_{k+1}$ intersects $B_i^i$ at finitely many points. But this implies that $B^i_{k+1}$ can not be singular along $C^i$. This argument applies to any irreducible component $B^k_j$ contracted along with $B^k_{k+1}$. Again Proposition 2.6 says that $B^k_{k+1}$ is normal and hence is a ${\Bbb P}^2$. A MET can be performed on $(X^k, B_{k+1}^k)$. After finitely many steps we obtain 
$$\Phi^m: X^m --\rightarrow X'$$ such that the proper transform $H^m$ of $H'$  is
numerically positive on $X^m$. This implies that $\Phi^m$ is an isomorphism
by Lemma 2.5.   

The proof of Theorem 1.1 is now complete.
\qed
\enddemo

The above proof has the following consequence on the uniqueness of METs.

\proclaim{Corrolary 3.2}
Let $\Phi: X --\rightarrow X'$ be a birational
tranformation of two hyperK\"ahler fourfolds which is
obtained by blowing-up a smooth center in each of $X, X'$. Then $\Phi$
is a Mukai elementary transformation up to  isomorphism.
\endproclaim

\heading
\S 4 Rational Hodge structure
\endheading

Theorem 1.1 also yields an isomorphism between  the rational
Hodge structures of two hyperK\"ahler fourfolds.

\proclaim{Corolary 4.1}Let $\Phi: X --\rightarrow X'$ be a birational morphism between two 
hyperK\"ahler fourfolds. Then
$\Phi$ induces an isomorphism between the rational Hodge structures of $X$ and $X'$.
\endproclaim

\demo{Proof}
Theorem 1.1 reduces the proof to the case when $\Phi: X --\rightarrow X'$ is a MET. Let $Z$ be the common blowdown of $X$ and $X'$ by collapsing the
 exceptional loci of $\Phi$ and $\Phi^{-1}$ to an isolated point $z$. 
Then the contraction $g: X \rightarrow Z$ ($g': X' \rightarrow Z$) is strictly semi-small
in the sense that $\text{dim} g^{-1}(z) = {1 \over 2} \text{codim} \{z \}$.
By the Beilinson-Bernstein-Deligne-Gabber decomposition theorem
(applied to the contraction $g$ and $g'$), we have the following
quasi-isomorphisms between perverse sheaves
$$ R g_* {\Bbb C}^\bullet_X \cong IC^\bullet (Z) \oplus IC^\bullet (z) [-4]$$
and 
$$ R g'_* {\Bbb C}^\bullet_{X'} \cong IC^\bullet (Z) \oplus IC^\bullet (z) [-4].$$
Both isomorphisms are compatible with rational Hodge decompositions by M. Saito's results. This proves the corollary.
\qed
\enddemo

\enddocument